# Coupled Routing and Charge Schedule Optimization of Electrified Delivery Truck Fleets: Feasibility Analyses

**Yared Tadesse Wendimagegnehu and Beshah Ayalew** Clemson University
**Andrej Ivanco** Allison Transmission Inc
**Habtamu Hailemichael** Clemson University

## Abstract

Electrifying truck fleets has the potential to improve energy efficiency and reduce carbon emissions from the freight transportation sector. However, the range limitations and substantial capital costs with current battery technologies imposes constraints that challenge the overall cost feasibility of electrifying fleets for logistics companies. In this paper, we investigate the coupled routing and charge scheduling optimization of a delivery fleet serving a large urban area as one approach to discovering feasible pathways. To this end, we first build an improved energy consumption model for a Class 7-8 electric and diesel truck using a data-driven approach of generating energy consumption data from detailed powertrain simulations on numerous drive cycles.

We then conduct several analyses on the impact of battery pack capacity, cost, and electricity prices on the amortized daily total cost of fleet electrification at different penetration levels, considering availability of fast charging at the depot. Findings indicate that at typical energy density of current battery technology, there is an optimal battery pack capacity that results from the contradicting effects of increasing pack capacity on cost, life span, weight and energy consumption. It is also observed that with currently improving trends in battery pack costs and availability of reduced electricity prices at the depot, such as with renewable microgrids, fleet electrification can become viable even at low levels of penetration.

## Introduction

Transportation electrification is set to play a significant role in reducing greenhouse (chiefly $CO_2$) emissions which contribute to global warming and climate change. Medium and heavy-duty trucks constitute 5% of the transportation sector but account for 24% of the transportation sector's $CO_2$ emissions [1]. In particular, logistics companies that deploy fleets of these vehicles in delivery and other services can benefit from the higher energy efficiency and therefore the lower operational cost of electric trucks. However, the high purchase costs of electric trucks, their limited driving range, and the requirement to transport heavy cargo over varying distances challenge their widespread adoption for heavy-duty applications [2], [3]. On the other hand, recent developments in high energy-density battery technologies hold the potential to make electric trucks technically and economically viable [2], [4], [5]. Furthermore, and in local delivery service, a single truck typically visits multiple customers, and each customer can be reached through different routes. When we also consider the charging needs for electric trucks because of their range limitations, the coupled optimization of route selections and charging schedules/durations (including partial recharges) presents an opportunity to identify optimal electrification levels and strategies. In this paper, we analyze the feasibility of electrifying a fleet of trucks for local delivery service considering variations in battery sizes and cost, as well as electricity price scenarios.

The optimal routing of electrified vehicles has been considered in several prior works, where it is sought to mitigate range anxiety through optimized route planning and charge scheduling [6], [7] [8], [9]. However, most works take a simplified approach to energy consumption modeling and focus on algorithmic aspects of rapidly solving for the optimal routes on standardized benchmark instances. However, economic and technical feasibility analysis of battery electric vehicles/trucks (BEVs) needs to consider realistic energy consumption models on real-world instances or instances closely aligned with actual service routes for these trucks.

There are indeed some pragmatic electric truck feasibility studies that analyzed the operational cost of BEVs, diesel trucks and trucks with other powertrains for inter-state, regional and local freight transportation services. Alonso et al. [3] studied the macro-scale techno-economic feasibility of alternative fuel trucks. In their work, the energy storage for both hydrogen and electric truck fleets was sized to match the effective energy of diesel trucks, and their operational costs for providing comparable services were analyzed. Their finding indicates a possibility of a higher adoption rate of electric trucks for local delivery services in Iceland [3]. Samet et al. [10] also conducted a comparative study of diesel and electric trucks based on levelized cost of driving. Their work concludes that electric trucks can be competitive for urban freight transport without any policy incentives. Similarly, analysis by NREL of a Manhattan beer distributors' truck duty cycle showed potential electrification benefits due to a shorter driving range, elongated idles hour, and lower driving speed in urban delivery services [11].

In a recent related work, Wang et al.in [4], [12] used a standard heterogenous fleet vehicle routing problem (VRP) formulation, with mixed-integer program (MIP) solvers, to conduct economic feasibility analysis of electrifying heavy-duty truck fleets. Their analysis related the feasibility of partial electrification of a delivery fleet to the diesel-to-electricity price ratio (USD/Gal/ USD/kWh). However, the fleet sizes considered were limited to under 20 customers, likely due to the computational difficulty of the exact MIP solution methods adopted, which are impractical for analysis of large-scale scenarios. They also assumed the availability of charging both at the depot and at customer nodes [4], the latter of which is less reflective of local delivery service scenarios at present.

The main contributions of this paper are: 1) an updated energy consumption model suited for routing optimizations, that considers explicit dependencies of energy consumption on road grade, travel speed between nodes, as well as varying weight, for both ICEV and BEV trucks. 2) analysis of the cost feasibility of electrified delivery fleets for large-scale scenarios, involving with many more customers than in previous works, by adopting the adaptive large neighborhood search (ANLS) metaheuristic algorithm to solve the coupled routing and charging problem [8], [13]. We also consider charging to be limited at the depot, where cheaper electricity price may be possible (such as with installation of microgrids integrating renewables and energy storage). The application of this study is to guide local delivery fleet operators and owners on some optimal strategies for electrifying their fleet.

## Energy Consumption Model

A good energy consumption model is crucial for estimating the energy cost of travel on any one link and tracking the state of energy (SOE) in routing decisions. A vehicle's energy consumption varies with road grade, traffic conditions, weather, driver behavior, and ambient temperature. Integrating all these considerations into the vehicle routing problem adds complexity and increases computational costs. Various simplification approaches are typically adopted to balance accuracy and computing cost. Some of these approaches include modeling travel energy as a multiple of driving distance with a constant coefficient [14], using energy consumption (ECC) coefficients based on the weight of the empty truck and freight weight [4], and using load and average speed-based energy consumption coefficient [9], [15], and splitting ECCs into acceleration-related and aerodynamics related coefficients [16].

In this work, we developed an energy consumption model from detailed (electrified and diesel) powertrain simulations. Several randomized combinations of urban and highway drive cycles were created with various average speeds and road grades. The average speed and road grade information for all routes in the service area is obtained from a map service provider [17]. Since the trucks transport freight that is delivered to different customers on the route, the vehicle's total weight along a route varies significantly. Similar to [4] [9], [15], we use empty truck and freight weight-related ECCs to consider these significant weight variations, but with the important difference that we explicitly model these ECCs' dependence on average speed and road grade. Numerous powertrain simulations are used to generate the empty truck ECC, $b_{ijk}$, and freight weight-related ECC, $a_{ijk}$, as a function of average speed and grade. We then used numerous drive cycle simulations to capture the effect of driving characteristics (kinetic intensity) on energy consumption at different road and traffic conditions. For instance, the urban drive cycle has a lower average speed and a high frequency of acceleration/ decelerations (high kinetic intensity), and highway drive has a high average speed and a lower frequency of acceleration and deceleration (low kinetic intensity) [18], [19], [20], [21]. The energy consumption $E_{ijk}$ of vehicle, $k$, while traversing link/edge $(i, j)$ of distance $d_{ij}$, transporting freight of weight, $w_{ijk}$ on that link, with an average speed of $v_{ij}$ is then modeled by:

$$E_{ijk} = a_k(v_{ij}, g_{ij}).w_{ijk}.d_{ij} + b_k(v_{ij}, g_{ij}).d_{ij} \qquad (1)$$

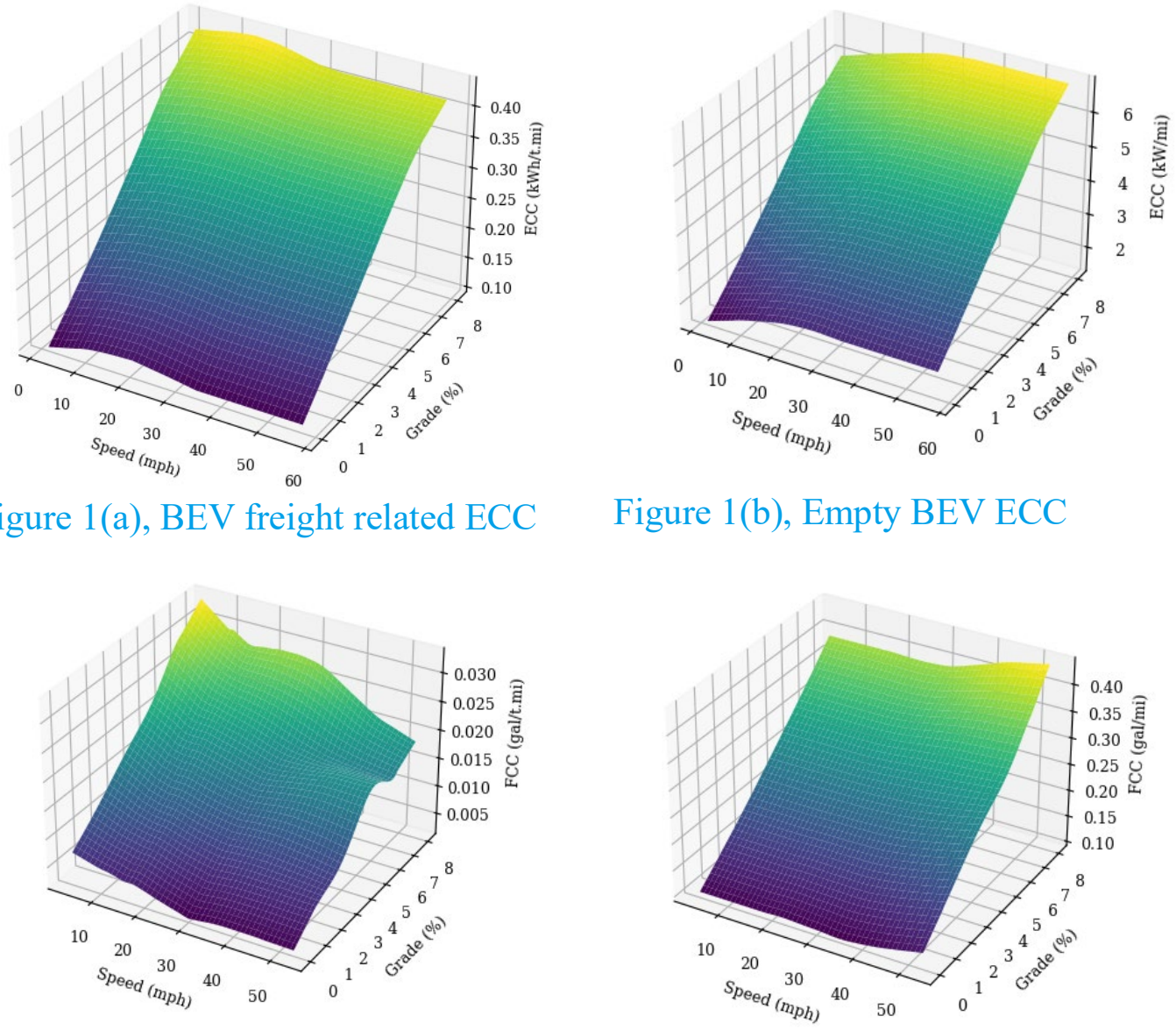

Figure 1(a), BEV freight related ECC

Figure 1(b), Empty BEV ECC

Figure 1(c), ICEV freight related FCC

Figure 1(d), Empty ICEV FCC

Figure 1 shows example ECC coefficient tables for an electrified truck model (top row) and the analogous Fuel Consumption Coefficient (FCC, in gal/mile or gal/ton.mile) for the diesel truck (bottom row) we used for our analysis. Note that as the battery sizes are increased, there is an increase in corresponding curb weight for the electric trucks and so the ECC coefficient tables $b_{ijk}$ are updated accordingly. When expressed in the same units, the ECC coefficients for the BEV are lower than that of the ICEV by as much as a factor of 2.5 or more.

Table 1 Definitions of variables

| $c_k^{ed}, c_k^{ec}$ | Recharging/refueling cost of truck-k (for BET depot and remote charging price are different but for diesel truck both are equal) |
|---|---|
| $c_d, c_v, c_m$ | Driver, levelized capital cost and maintenance of cost of a truck per mile |
| $E_{ijkr}$ | Energy consumption of truck-k while traversing edge $(i, j)$ on its $r^{th}$ trip |
| $E_{ikr}, E_k^{min}, E_k^{max}, E_k^{cap}$ | Arrival energy, lower limit, upper limit, and energy capacity of truck-k, respectively |
| $x_{ijkr}$ | Binary decision variable if truck-k traverse edge-$(i, j)$ on its r-trip or not |
| $u_{kr}$ | Binary indicator variable if arrival state of energy is below minimum or not |
| $P_{xk}$ | Auxiliary system energy or fuel consumption rate |
| $d_{ij}, t_{ij}, v_{ij}, g_{ij}$ | Link-$(i, j)$ fastest route distance, travel time, average speed, and average grade |

| $\tau_{aikr}, \tau_{sikr}$, $t_{si}, t_{ai}, t_{bi}$ | Arrival time, service start time, service duration, lower and upper bound of time window at node-$i$ |
|---|---|
| $q_i, w_{ijkr}, Q_k$ | Customer demand, freight weight, and payload capacity in tons |

## Problem Formulation

The optimization problem is formulated as a heterogeneous vehicle routing problem. The VRP is composed of two kinds of trucks: electric trucks and diesel trucks. Electric trucks are allowed to partially charge at the depot, and diesel trucks are allowed to fully refuel at gas stations (with no significant time). Each truck is allowed to be dispatched multiple times to meet customer demand. Each customer's nodes have a loosely defined service time window and are characterized by demand, arrival time and service time. The aim is to minimize the total cost of serving all the customers' demands. See details on the customer demand modeling in the Instance Generation section below.

### *Objective function*

The objective function comprises daily operational costs related to energy, driver, capital, and maintenance costs:

$$J = \sum_{r\in R}\sum_{k\in K}\sum_{(i,j)\in A} c_k^{ed} E_{ijkr}.x_{ijkr} + \sum_{r\in R}\sum_{k\in K}\sum_{i\in N_s} E_{cikr}(c_k^{ec} - c_k^{ed})x_{ijkr} + \sum_{r\in R}\sum_{k\in K}\sum_{(i,j)\in A} c_k^{ed} P_{xk} t_{ij}.x_{ijkr} + \sum_{r\in R}\sum_{k\in K}\sum_{(i,j)\in A} c_d d_{ij} x_{ijkr} + \sum_{r\in R}\sum_{k\in K}\sum_{(i,j)\in A} (c_v + c_m) d_{ij} x_{ijkr} \quad (2)$$

Definitions of all notations is given in Table-1. In equation (2), the first three terms are travelling energy, differential recharging /refueling cost, and auxiliary energy costs. The fourth term of the objective function is driver cost, and the last term is capital expense and maintenance cost, all of which are amortized per mile of travel, considering the service life of the vehicles. Optimization of the above objective function is subjected to the following network flow constraints.

$$\sum_r\sum_k\sum_j x_{ijkr} \geq 1, \quad \forall r\in R, \forall\, k\in K, \forall (i,j)\in A, \&\, i\neq j \quad (3)$$

$$\sum_i x_{ijkr} \leq 1, \quad \forall r\in R, \forall\, k\in K, \forall (i,j)\in A, \&\, i\neq j \quad (4)$$

$$\sum_i x_{ijkr} - \sum_j x_{ijkr} = 0\ , \forall r\in R, \forall k\in K, \forall (i,j)\in A, i\neq j \quad (5)$$

$$\sum_r\sum_k\sum_i x_{ijkr} \leq n_{rv} \quad \forall r\in R, \forall k\in K_E, \forall i\in N_{cs} \quad (6)$$

$$\sum_r\sum_k\sum_i x_{ijkr} \leq n_{rv}\ , \forall r\in R, \forall k\in K_D, (i,j)\in A, \forall i\in N_{gs} \quad (7)$$

$$\sum_i x_{ijkr} \geq u_{kr}, \quad \forall r\in R, \forall k\in K, \forall (i,j)\in A, \forall j\in \{N_{cs}\cup N_{gs}\} \quad (8)$$

$$E_{nkr} \leq E_k^{min} - \epsilon + E_k^{cap}(1-u_{kr})\ , \forall r\in R, \forall k\in K, \forall u\in\{0,1\} \quad (9)$$

$$E_{nkr} \geq E_k^{min} + \epsilon + E_k^{cap} u_{kr}, \forall r\in R, \forall k\in K, \forall u\in\{0,1\} \quad (10)$$

$$E_{jkr} \leq E_{ikr}. + P_i\tau_{cikr} - E_{ijkr}x_{ijkr} + (1-x_{ijkr})0.8E_{cap}, \ \forall r\in R, \forall k\in K, \forall (i,j)\in A \quad (11)$$

$$t_{ai} \leq \tau_{sikr} \leq t_{bi} \quad i:(i,j)\in A \quad (12)$$

$$\tau_{sikr} + \tau_{cikr} + t_{si} + t_{ij} - \tau_{ajkr} \leq (1-x_{ijkr}).M_{ij}\ \forall r\in R, \forall k\in K, \forall (i,j) \quad (13)$$

$$Q_k \geq \sum_i q_i\ , \forall r\in R, \forall\, k\in K, \forall (i,j) \quad (14)$$

$$\sum_i w_{ijkr} = q_j + \sum_j w_{ijkr}\ , \forall r\in R, \forall\, k\in K, \forall (i,j) \quad (15)$$

The first constraint (3) ensures that all customers are served, but not all customers are served on the first trip as specified in constraint (4). Constraint (5) is a flow conservation constraint that enforces all arcs entering a node shall leave it. Constraints (6) and (7) limit the number of replenishing (charging, refueling) station visits to electric and diesel trucks, respectively. Each truck can visit the charging/refueling station at most $n_{rv}$ times in a single trip. Constraint (7) ensures that the vehicle visits at least one refueling/charging station if its arrival energy is below the minimum. This constraint can be adjusted to set energy reserves in mixed integer programming formulation of vehicle routing. The indicator variable $u_{kr}$ is a big-M constraint in equation (9-10). Constraint (11) tracks energy remaining energy of the truck at nodes. Constraint (12) ensures the customers are visited in the predefined time window. Constraint (13) tracks the vehicle arrival and departure time at a node. Constraint (14) limits the vehicle freight load within the truck's payload capacity. Constraint (15) tracks the total load while traversing the edge $(i,j)$, and it eliminates subtours.

The above mixed-fleet and multi-trip vehicle routing problem with charging and time windows is indeed a mixed-integer programming problem. However, this formulation is very difficult to solve in reasonable time with state-of-the art exact methods for MIPs when we consider realistic problem sizes with several 10s or hundreds of customer nodes (large instances). In this work, we modified the Adaptive Large Neighborhood Search (ANLS) algorithm and used it to solve the coupled routing and charge scheduling problem for electrified fleets serving large instances.

### *Adaptive Large Neighborhood Search (ALNS)*

ALNS, a comprehensive large neighborhood search meta-heuristic algorithms, is commonly used for efficiently solving large-scale capacitated vehicle routing problems[8]. The algorithm begins with a feasible initial solution and selectively applies destroy (Removal) and repair (Construction) operators [22], [23] to generate a better neighborhood solution [22], [23]. The new solution is compared with the current solution and accepted or rejected based on predefined acceptance criteria. Satani et al. have studied the performance of the ALNS solution acceptance method for various applications using the WIXCON test [22]. They found a Linear Record-to-record travel (RRT) acceptance method with fixed-end outperforms the simulated annealing acceptance method for capacitated vehicle routing problems [22]. Thus, in our ALNS implementation the linear RRT is used for acceptance, and Roulette wheel algorithm for operator selection and ranking. The Removal and Construction heuristics are briefly described next.

### *Removal heuristics*

The removal heuristic of the ALNS algorithm eliminates customer nodes from the constructed routes. These operators leave a partially destroyed route to repair and build better routes for the construction heuristics. We have implemented location-based removal, random removal, worst removal, route removal and Shaw removal heuristics. The worst removal is extended to the worst distance, worst cost, and worst time removal operator. Similarly, Shaw removal that removes related nodes from routes is extended to time, space, truck, and demand Shaw removal [22], [24], [25].

### *Construction heuristics*

Construction heuristics build a better route by inserting removed nodes into a partially destroyed route. The repair heuristics can be a sequential heuristic, which constructs one route at a time; parallel insertion heuristics develop multiple routes simultaneously [22].The ALNS implementation includes greedy repair with and without noise, random insertions, and regret-2 repair heuristics.

## Instance Generation

We consider delivery service scenarios derived from the grocery store customer locations and electricity pricing for the greater Atlanta, Georgia Area (Fig. 2a). Once the locations of the depot and charge/refuel stations are fixed, the customer locations and demand are generated randomly within this potential service area (Fig. 2b). The randomization aims to capture diverse route characteristics without being limited to a fixed set of customers in our feasibility analysis. To this end, we consider a normal distribution of customer geographic locations centered at the depot and with standard deviation of one third of the width of the service area box, as measured in degrees of latitude and longitude. Figure 2 an example disposition of randomly generated customer locations (b), motivated by the actual potential customer locations from actual grocery store locations (a). In a similar vein, we randomized customer demand (delivery load) by uniformly sampling it b/n 1 and 3.5 tons. An instance can be considered as a daily variation considering both the geographic distribution of the customers and their delivery demands. The randomization varies both aspects between instances.

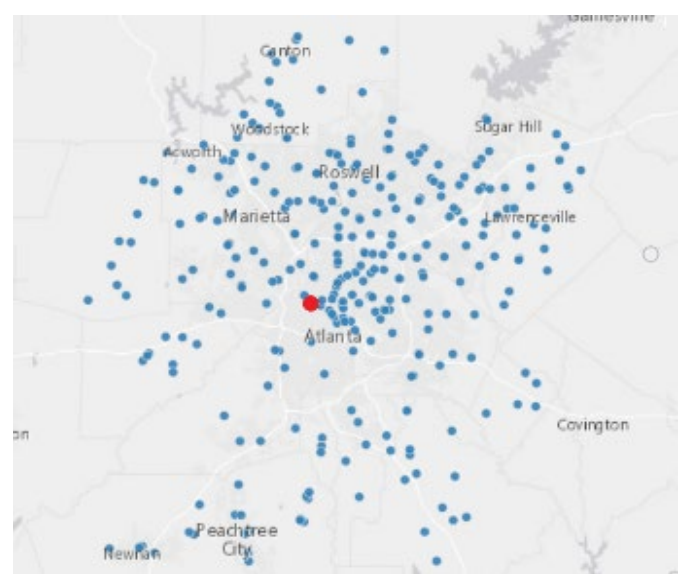

Figure 2(a) Actual grocery store customers in Greater Atlanta Area [26]

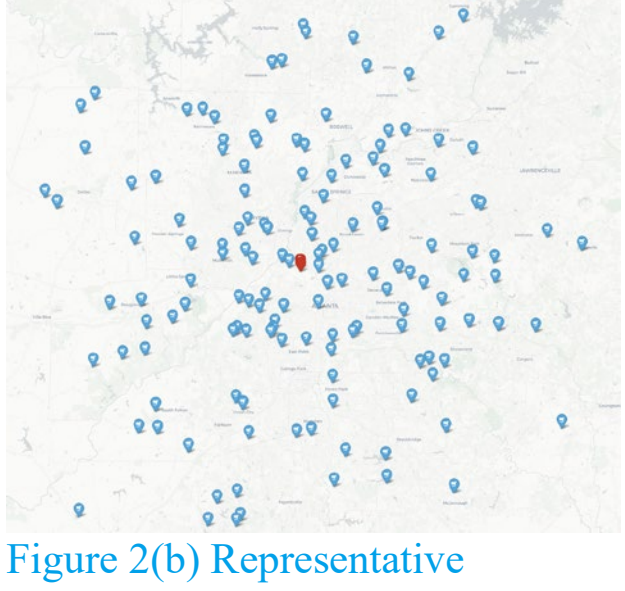
Figure 2(b) Representative randomized instances.

Table 2 Vehicle Parameters

| Vehicle | Payload [tons] | Values |
|---|---|---|
| Diesel truck | Payload capacity[tons] | 23.8 |
| | Empty truck weight[tons] | 13.21 |
| | Fuel tank [gal] | 240 |
| Electric truck | Payload capacity[tons] | 22.11 |
| | Empty truck weight[tons] | 14.18 |
| | Battery size [kWh] | 440 |

Table 3. Optimization Parameter Settings

| Parameters | Values |
|---|---|
| Depot default electric price [USD/kWh] | 0.115[27] |
| Diesel fuel price (avg) [USD/gal] | 3.47-4.11 |
| BEV maintenance [USD/mi] | 0.141 [28] |
| ICEV maintenance [USD/mi] | 0.196 [28] |
| BEV truck cost [1000 USD] | 248 [28] |
| ICEV truck cost [1000 USD] | 129 [28] |
| BEV Li-ion battery life [k.mi] | 300 for 440kWh [29], [30] |
| Truck lifespan [kmi] | 750 [30] |
| Diesel on-board $CO_2$ emission [kg/gal] | 750 [31] |

## Results and Discussion

In our first set of analyses below, we set the ALNS algorithm termination criteria to 10,000 iterations and ran it on ten randomized instances of size 150 customers, executing it five times on each instance. Each instance has a total load demand of 400-450 tons. We considered a fleet size of 10 while analyzing costs listed in the objective function with varying electricity prices, battery costs, and the percentage penetration of electric trucks as we executed the ALNS algorithm for the coupled optimization of routing and charging of the mixed fleet.

### *Effect of Electricity Pricing and Battery Cost*

Electric trucks have a potential advantage over diesel trucks as having a more efficient powertrain (lower ECC), lower maintenance cost and zero onboard carbon emissions, but they also cost more. We analyzed the total daily amortized cost of electrified fleets at various electricity prices (current (late 2024) prices in Georgia being between \$0.070/kWh and \$0.115/kWh, and even lower prices are assumed feasible with renewables at the depot). In this study, truck charging infrastructure investment is not considered. The analysis is also extended to include the effect of battery pack price, where we considered two battery pack prices reflective of current market trends [32]. The results of the coupled routing and charging optimization are then summarized with the detailed daily cost plots at different levels of electric truck percentage penetration in the fleet. Both the total amortized daily (averaged per instance) cost and the variations are shown. In all cases, we will see that there is only a small variability between instances in both the objective function components and the overall daily cost, showing the repeatability of the analysis results from ANLS.

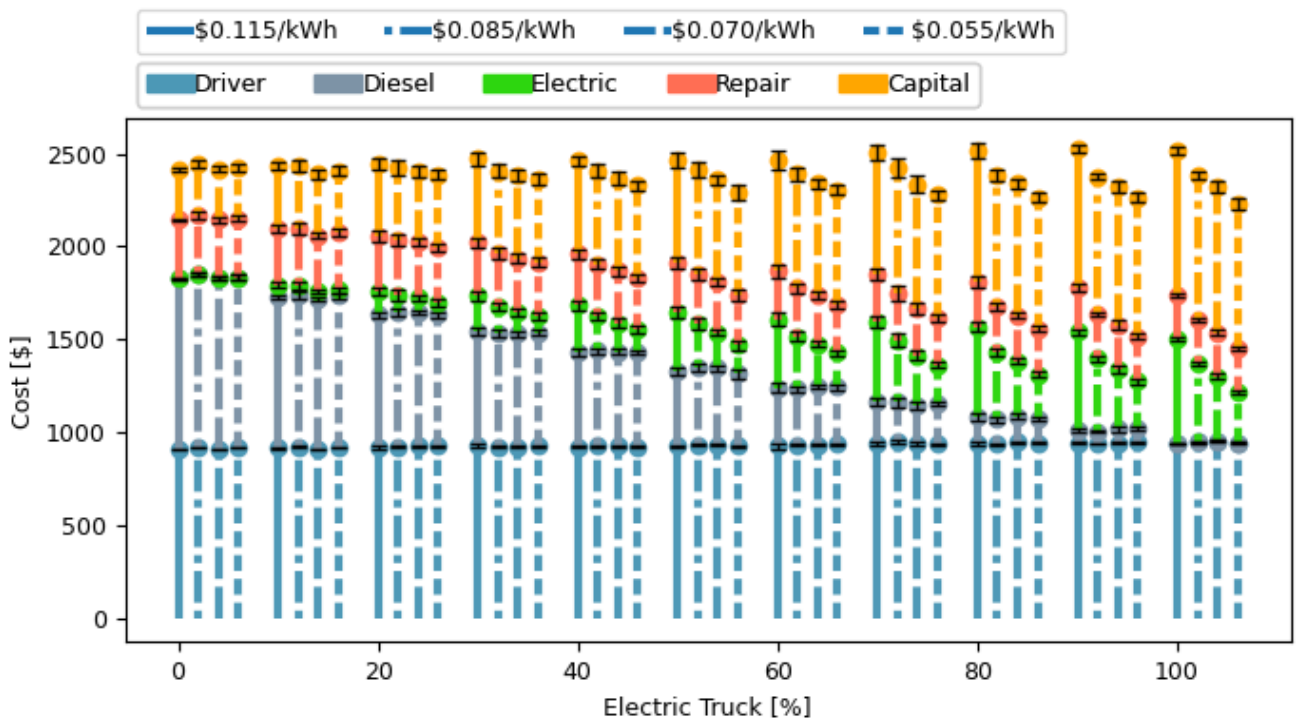

Figure 3. Average daily cost of serving 150 customers with 10 trucks at a battery pack price of \$240/kWh and fuel price of \$3.86/gal. Error bars show variation of the ALNS's solution between instances.

Figure 3 shows that at a battery pack price of \$240/kWh, the total daily cost increases with an increase in percentage penetration of electric trucks when the electricity price is above \$0.055/kWh, since the amortized capital costs are not able to be fully offset with electrification. It is only at the low electricity price of \$0.055/kWh that electrification becomes fairly competitive with the all diesel fleet (0% penetration), starting with around a 40% penetration. As already stated, this low price is likely possible with installation of microgrids (renewables and storage) at the depot. However, Figure 4 shows that for a battery pack price of \$153/kWh (which is currently close to the cheapest possibility listed in the market at \$139/kWh [32]), the benefits of fleet electrification in reducing the daily cost are possible at around a \$0.070/kWh electricity price and starting even with a low 10% penetration.

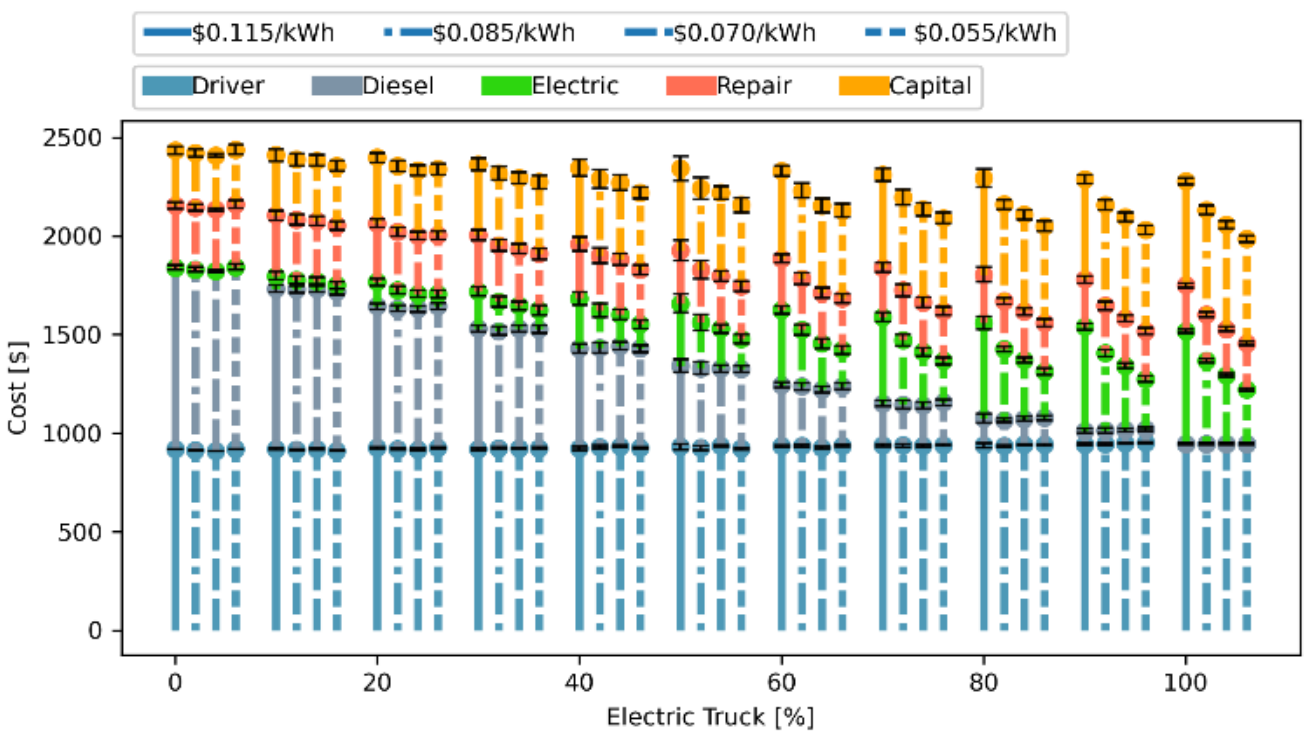

Figure 4. Average cost of serving 150 customers with 10 trucks at battery size price of \$153/kWh and fuel price of \$3.86/gal.

Figure 5 shows the accompanying daily mileage distribution. As could be expected from the range limitation of electric trucks, more total mileage is involved with increasing electric truck penetration and more % electric mileage is involved than the linear expectation by the fleet composition. For example, with a 60% penetration, close to 70% of the mileage is electric. However, although some of this could be increased electric mileage from charging trips to the depot; and in our formulation trucks pick up additional demand on these trips as well and so are productive mileage.

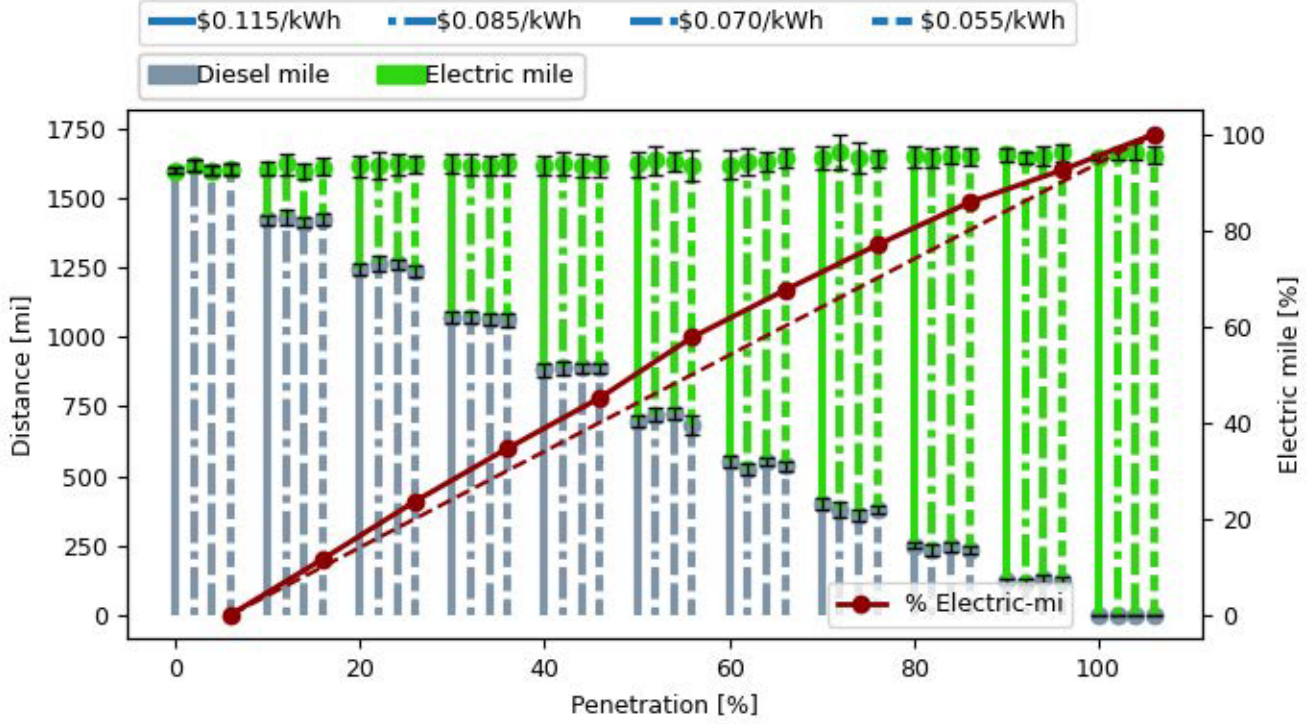

Figure 5. Total daily mileage of the fleet to serve the 150 customers, battery pack price of \$240/kWh and fuel price of \$3.86/gal. The electric mile % is shown only for a price of \$0.055/kWh (compare with the dotted line showing the linear split).

Another customary way to look at the above results is to consider the cost per mile, which is summarized in Figure 6. In general, except when the electricity price is $\$0.115/kWh$, the operational cost per mile drops with increasing percentage penetration. This is due to the general decreases in total cost and accompanied increases in mileage observed for the different electricity price points. Viewed in this way, if all the electric mileage is considered productive, the cost per mile may be competitive even at \$0.085/kWh at even 20% penetration.

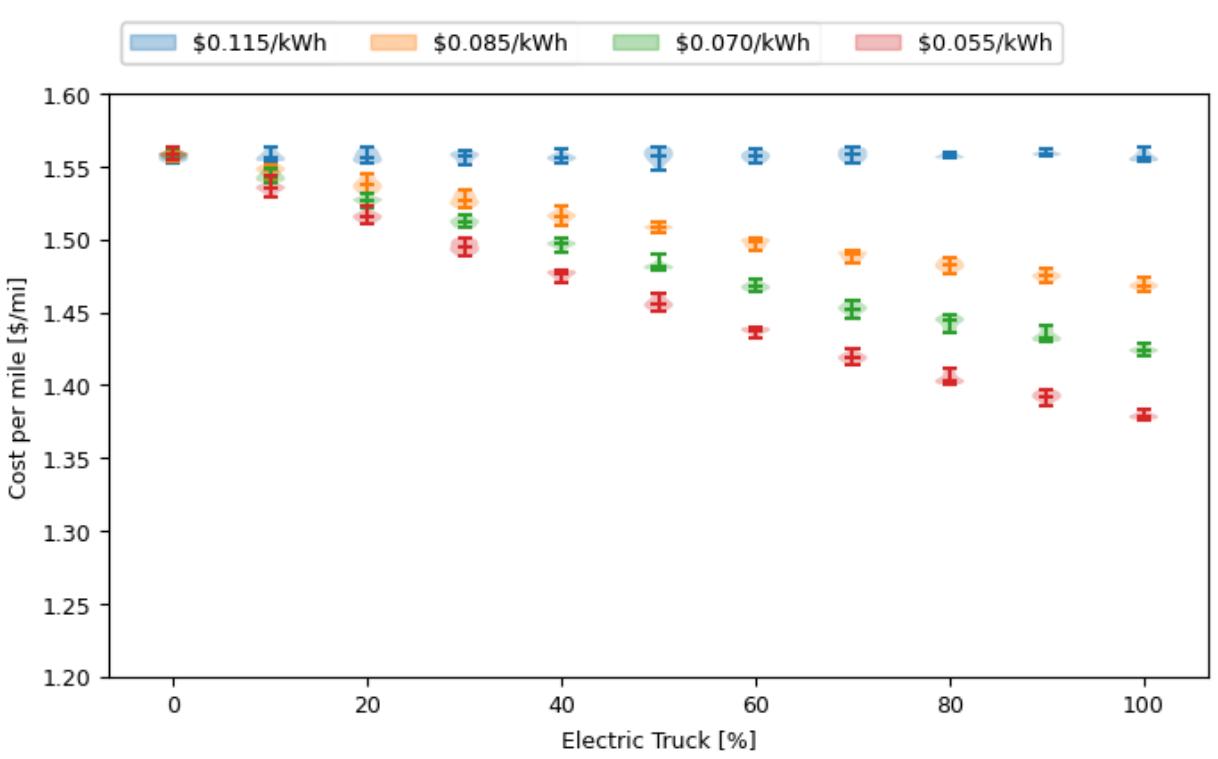

Figure 6. Operational cost per mile of serving a 150 customer with 10 trucks at battery pack price of \$240/kWh and fuel price of \$3.86/gal.

Overall, these results show that making delivery fleet electrification competitive requires lower electricity prices to even take advantage of lower battery pack prices, which are trending in the right direction in the marketplace [32]. An equivalent interpretation is that coupled routing and charge scheduling optimization for a delivery fleet with just depot charging could make electrification competitive if diesel prices increase substantially relative to the electricity prices available at the depot.

### *Effect of Battery Capacity*

Large battery capacity or size gives longer driving range and possibly longer life mileage, $L_{bat}$. However, it comes with increased truck curb weight, reduced payload capacity, $Q_{bev}$ (for the same GVW), and increased purchasing cost of the vehicle, $C_{bev}$. In these specific experiments, the battery cost is set to $\$240/kWh$, and the base truck battery is 440kWh.

I. Purchasing cost: We assume a linear interpolation estimate of the capital cost as a function of battery pack cost.

$$C_{BEV}(E_{new}) = C_{base} + 240(E_{new} - E_{base}) \tag{16}$$

II. Payload capacity: For a given class of truck, a new payload capacity, $Q_{BEV}$, with a new battery size, $E_{new} \geq E_{base}$ is updated based on the energy density of the lithium-ion battery, $\rho_{bat}$. The energy density of modern lithium batteries ranges from 260-280 kWh/ton [33], and we assumed 270 kWh/ton.

$$Q_{BEV}(E_{new}) = Q_{base} - (E_{new} - E_{base})/\rho_{bat} \tag{17}$$

III. Empty truck energy consumption coefficient ($b^*_{ijk}$). This is updated and re-generated for the energy consumption model.

$$b^*_{ijk} = b_{ijk}(v_{ij}, g_{ij}) + a_{ijk}(v_{ij}, g_{ij}).(E_{cap} - E_{base})/\rho_{bat} \tag{18}$$

IV. Battery lifespan in miles

A larger battery allows a longer driving range and decreases the charging frequency, which can improve battery life. For this specific case, we linearly interpolate the BEV powertrain warranty provided by eCascadia to estimate a truck battery life, $L_{bat}$, with increased battery size, $E_{cap}$, in kWh [29].

$$L_{bat}(E_{cap}) = 1.0204(E_{cap} - 291) + 150\ kmi \tag{19}$$

The capital expense of an electric truck is computed assuming battery pack replacements to provide a service lifespan for the electric truck that is equivalent to that of a diesel truck. The diesel truck can last up to 750kmi without an engine overhaul [1], [31].

In the following analyses of the effect of battery pack capacity/size, we consider an instance size of 79 customers with a total demand of 77-114 tons served with a fleet of 8 trucks. Figures 8 and 9 summarize the results on the effect of battery capacity on the average daily total cost and other cost components. While larger battery capacities than the current default of 440kWh seem to reduce the total daily cost, there is an optimal battery capacity of around 528kWh for which the fleet's cost per mile is minimal from the capacity/sizes considered. This is due to the trade-off in the effects of increasing battery capacity as described above.

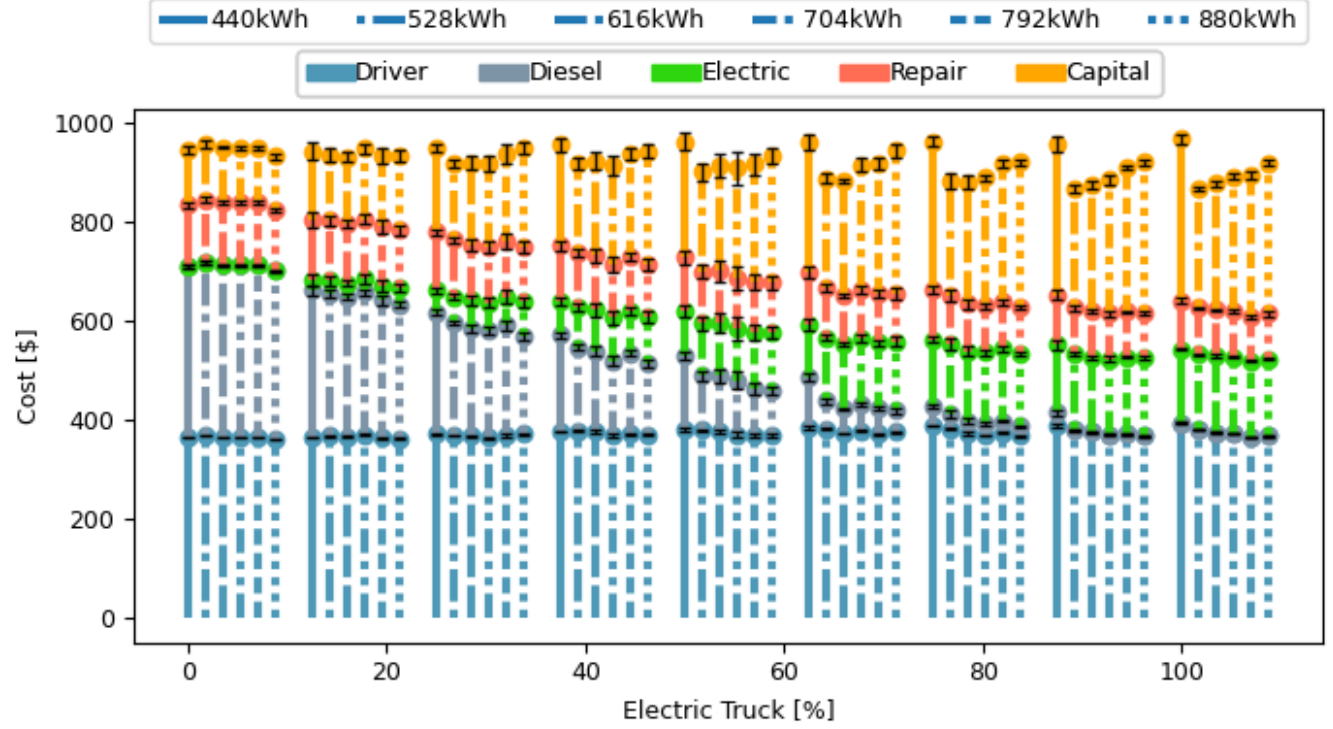

Figure. 8 Average cost per instance of serving 79 customers with 8 trucks and electricity price $0.075/kWh (battery pack price of $240/kWh, fuel price of $3.65/gal and electricity price of $0.075/kWh).

Figure 9 illustrates the distribution of the average cost per mile for 8 instances. The average cost per mile captures the effect of the total driving distance of the fleets associated with the percentage penetration of electric trucks as mentioned above. Figure 10 shows that increasing electric miles with increasing penetration (which still holds even if more customers are being served with fewer recharging trips). This explains the decreasing trend in the cost per mile for all battery pack capacities considered. We observe from Figure 9 that there is an optimal battery capacity of around 528kWh for these instances when measured with the cost-per mile metric as well. The result in Figure 10 also shows, for this particular scenario, one truck does not even need to be dispatched for the fleet with 87.5% penetration for battery sizes greater or equal to 528kWh.

Finally, it is important to note that the optimal capacity will be different for different demand distributions (customer sizes, geographic and cargo demand distributions).

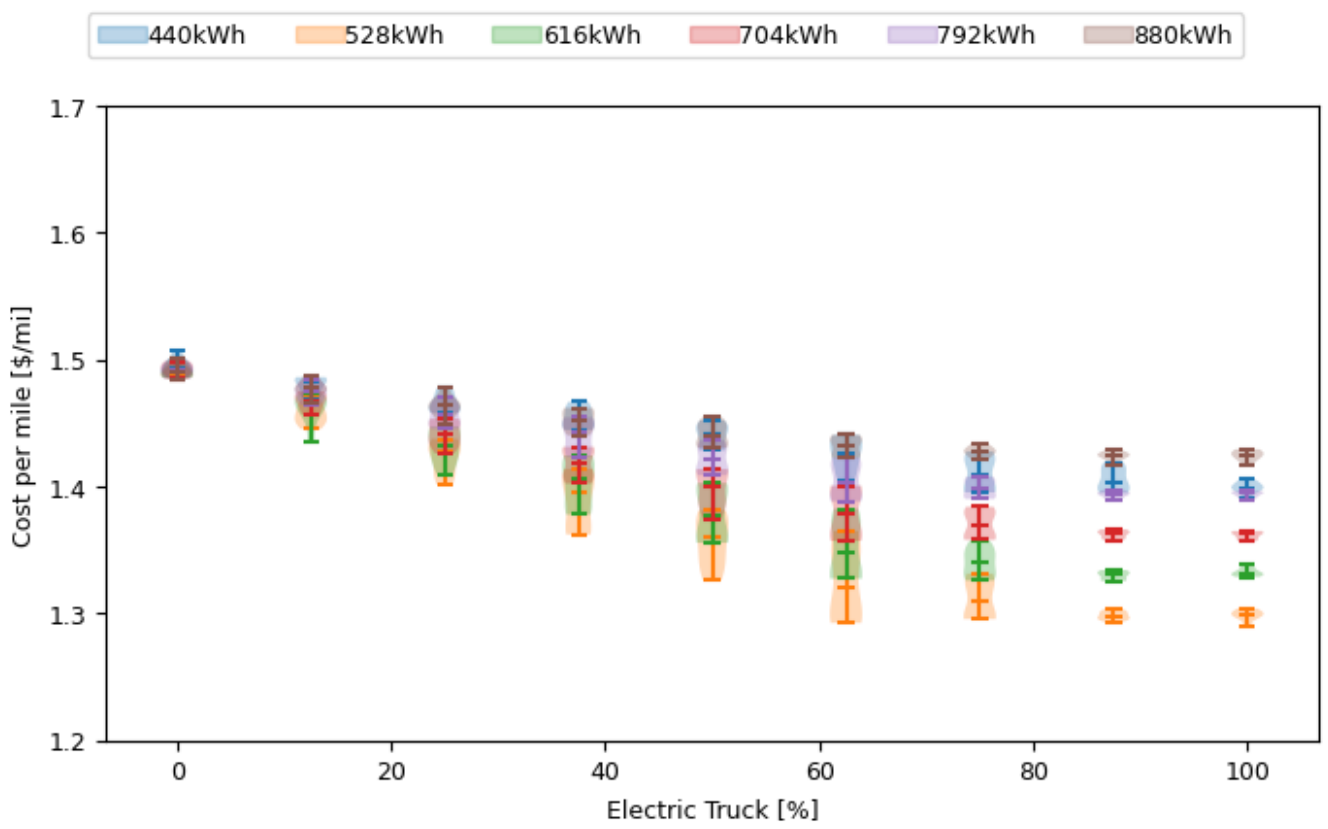

Figure 9. Cost per mile of serving 79 customers with 8 trucks (battery pack price of $240/kWh, fuel price of $3.65/gal and electricity price of $0.075/kWh).

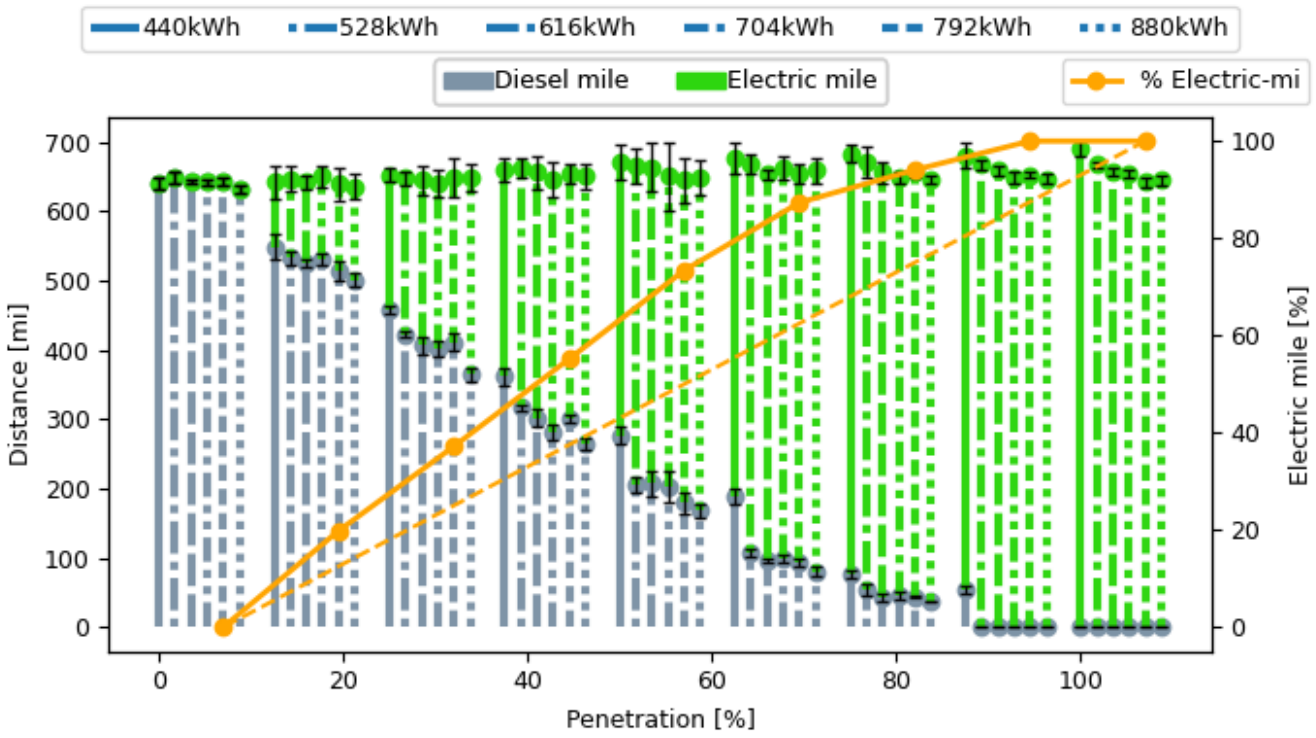

Figure 10. Total daily mileage of the fleet to serve the 79 customers (battery pack price of $240/kWh, fuel price of $3.65/gal and electricity price of $0.075/kWh).

## Summary/Conclusions

The electrification of local delivery trucks has a clear potential for improving the energy efficiency and climate impact of freight transportation. This paper analyzes the total amortized daily cost of using electrified fleets for large-scale delivery services via a coupled routing and charge scheduling optimization approach. Using adaptations of the ANLS algorithm, with detailed diesel and electric truck energy consumption model, and considering large-scale delivery scenarios served with heterogenous fleets, we conducted analysis on the effects of electricity price, battery pack purchasing price and battery capacity. The results show that electric truck competitiveness can be increased by optimally sizing the battery pack and that this can lead to cost parity or improvement compared to an all diesel fleet. We also noted that partial electrification can become competitive at lower differentiated electricity prices possible at the depot and is becoming even more viable with the decreasing trends in battery pack prices. These potential benefits are observed without explicit considerations of policy incentives or subsidies that favor electrification. The practical implication of the findings is logistic companies and fleet owners and operators can realistically pursue cost competitive electrification strategies, by optimally sizing the battery packs and taking steps (such as with renewables, or negotiating with utilities) for reduced electricity pricing for depot-based charging of their electrified fleets.

The analysis has some limitations: It did not consider potential queuing issues at the depot when the number of charging ports is limited, and it did not consider the cost of charging infrastructure. The impacts of thermal/climatic considerations is another aspect that needs further study. Furthermore, although the formulation includes remote or enroute charging possibilities, the analysis did not include those scenarios. Furthermore, vehicle and fleet sizing issues are not covered. In addition, the purchasing costs of electric and diesel trucks in Table-2 are based on the referenced literature [28]; actual prices in the market could differ substantially from those listed due to localized incentives and other market factors. Finally, differences in residual costs of BEVs vs. ICEVs beyond the 750k mile service life are not explicitly analyzed. These aspects will be addressed as future work.

## Acknowledgments

The authors acknowledge the support provided by Allison Transmission Inc for a research project on this topic. The results and opinions expressed herein are of the authors only.

## Definitions/Abbreviations

| | |
|---|---|
| **ALNS** | Adaptive Large Neighborhood algorithm |
| **BEV** | Battery Electric Vehicle |
| **CapEx** | Capital expense |
| **ECC** | Energy consumption coefficient. |
| **ECR** | Energy consumption rate |
| **FCC** | Fuel consumption coefficient |
| **ICEV** | Internal combustion engine vehicle |
| **LCOD** | Levelized cost of driving |
| **PPET** | Percentage penetration electric truck |
| **VRP** | Vehicle routing problem |